\documentclass[11pt, reqno]{amsart}

\usepackage{amsmath,amssymb,amsthm, epsfig}
\usepackage{amsmath}
\usepackage{amssymb}
\usepackage{color}
\usepackage{mathtools}
\usepackage{bbm}
\usepackage{calrsfs}
\DeclareMathAlphabet{\pazocal}{OMS}{zplm}{m}{n}
\usepackage{epsfig}
\usepackage[mathscr]{eucal}
\usepackage[latin1]{inputenc}
\usepackage{enumerate}
\usepackage[shortlabels]{enumitem}

\usepackage{xcolor}
\usepackage{hyperref}
\usepackage[nameinlink]{cleveref}
\hypersetup{
	colorlinks,
	linkcolor={blue!50!black},
	citecolor={blue!50!black},
	urlcolor={blue!50!black}
}

\usepackage{blindtext}

\newtheorem{theorem}{Theorem}
\newtheorem{definition}{Definition}
\newtheorem{lemma}{Lemma}
\newtheorem{proposition}{Proposition}

\numberwithin{equation}{section}
\numberwithin{theorem}{section}
\numberwithin{lemma}{section}
\numberwithin{corollary}{section}
\numberwithin{remark}{section} \numberwithin{proposition}{section}
\numberwithin{definition}{section}
\numberwithin{assumption}{section}

\usepackage[margin=1in]{geometry}

\newcommand{\dd}{\mathrm{d}}
\newcommand{\supp}{\operatorname{supp}}

\newcommand{\loc}{\operatorname{loc}}

\newcommand{\Div}{\operatorname{div}}

\newcommand{\bv}{\operatorname{BV}}
\newcommand{\tr}{\operatorname{Tr}}
\newcommand{\bl}{\operatorname{LD}}
\newcommand{\sym}{\operatorname{Sym}}
\newcommand{\pseudo}[1]{{\left\vert\kern-0.25ex\left\vert\kern-0.25ex\left\vert #1 
		\right\vert\kern-0.25ex\right\vert\kern-0.25ex\right\vert}}

\usepackage[foot]{amsaddr}

\begin{document}

\title[Regular Lagrangian flows for Hörmander singular kernels and \texorpdfstring{$\bl$}{} vector fields]{Quantitative estimates for flows of regular Lagrangian flows for Hörmander singular kernels and \texorpdfstring{$\bl$}{} vector fields}

\author[H. Borrin]{Henrique Borrin$^\dagger$}
\address{$^\dagger$Faculdade de Filosofia, Ci\^{e}ncias e Letras de Ribeir\~{a}o Preto, USP-Universidade de São Paulo \\ Address: Avenida Bandeirantes, 3900. Ribeir\~{a}o Preto, SP, Brazil. Zip Code 14040-901. ORCID: 0000-0003-4670-4444.}
\email{henriqueborrin@usp.br (Corresponding author)}

\thispagestyle{empty}

\begin{abstract}
In this paper, we obtain quantitative estimates of regular Lagrangian flows associated to vector fields whose derivative can be written as convolution of a fundamental singular kernel satisfying the ``Hörmander'' condition convoluted with summable function in spacetime.

\bigskip

\bigskip

\end{abstract}

\maketitle

\section{Introduction}
In this paper, we shall establish the well-posedness of solutions for the flow equation
\begin{equation}\label{flow}
	\begin{cases}
		\displaystyle\frac{\dd}{\dd t}\boldsymbol{X}(t,s,x)=\boldsymbol{b}(t,\boldsymbol{X}(t,s,x))\quad &\text{for } (t,x)\in (0,T)\times\mathbb{R}^d;\\
		\boldsymbol{X}(s,s,x)=x \quad&\text{for } x\in\mathbb{R}^d
	\end{cases}
\end{equation}
for vector fields $\boldsymbol{b}:[0,T]\times\mathbb{R}^d\rightarrow \mathbb{R}^d$ with the ever-present growth assumption introduced in \cite{dipernalions}
\begin{equation}\label{growth}
	\begin{split}
		\frac{\boldsymbol{b}(t,x)}{1+|x|}=\boldsymbol{b}^1(t,x)+\boldsymbol{b}^2(t,x),\text{where } \boldsymbol{b}^1\in L^1([0,T]\times\mathbb{R}^d),\, \boldsymbol{b}^2\in L^1([0,T];L^\infty(\mathbb{R}^d)),
	\end{split}
\end{equation}
and whose space symmetric derivative can written as a convolution:
\begin{equation}\label{derivativeconvolution}
	\partial_j\boldsymbol{b}^i(t,x)+\partial_i\boldsymbol{b}^j(t,x)=\sum^m_{k=1}\Gamma_{ijk,R}\ast g^{ijk,R}(t,x)\quad \text{for } (t,x)\in(0,T)\times B_R,
\end{equation}
where $g^{ijk,R}\in L^1([0,T]\times\mathbb{R}^d)$ for each index and $\Gamma_{ijk,R}$ is singular kernel \textit{à la} Calderón-Zygmund theory for each index, in the sense summarized in \cite[Definition 2.7, Proposition 2.12]{bouchutcrippa}.
\begin{definition}[Fundamental singular kernel]\label{singularkernel}
	\textnormal{We say that $\Gamma\in L^1_{\loc}(\mathbb{R}^d\setminus \{0\})$ with bounded Fourier transform, i.e. $\mathcal{F}(\Gamma)\in L^\infty(\mathbb{R}^d)$, is a fundamental singular kernel if it holds
	\begin{equation}\label{hormander}
	C_0\coloneqq\sup_{z\neq 0}\int_{|y|>2|z|}|\Gamma(y-z)-\Gamma(y)|\,\dd y<\infty.
\end{equation}
Equivalently, a fundamental singular kernel $\Gamma$ is a function satisfying \eqref{hormander} and 
\begin{equation}\label{cancellation}
	C_1\coloneqq\sup_{r>0}\int_{r<|y|<2r}|\Gamma(y)|\,\dd y \quad\text{and}\quad  \left|\int_{r<|y|<R}\Gamma(y)\,\dd y\right|<\infty
\end{equation}
for all $0<r<R<\infty$ for some constant $C>0$.}
\end{definition}

Using the same argument as in the comment after main theorem in \cite{nguyen}, one may consider vector fields in $L^1((0,T);\bl_{\loc}(\mathbb{R}^d))$, where we write as in \cite{conferenceambrosio} the space
\[\bl(\Omega)=\left\{\boldsymbol{u}\in L^1(\Omega;\mathbb{R}^d): \partial_j \boldsymbol{u}^i+\partial_i \boldsymbol{u}^j\in L^1(\Omega)\text{ for each } i,\,j \in 1,\dots,d\right\}.\]

The bibliography comprehending the well-posedness of \eqref{flow} for non-smooth vector fields is quite extensive, with major contributions (but not limited to) in  \cite{ambrosio,inftydim,inftydim2,anisotropic,borrinvlasovmaxwell,borrinmarcon,bouchutcrippa,symmetricderivative,crippadelellis,nussenzveig,partiallyregular,pointcharge,dipernalions,specularxavier,inftydim3,miotsingularset,nguyen}. We remark that \cite{symmetricderivative} is the main inspiration of this article, for it generalizes the technique pioneered in \cite{dipernalions}, which was limited to $L^1((0,T);W^{1,1}_{\loc}(\mathbb{R}^d))$ vector fields, for the more general class $L^1((0,T);\bl_{\loc}(\mathbb{R}^d))$ by a simple assumption of restricting the mollifier to radial functions. However, the Lagrangian variation for such vector fields was not available; we shall specify our contributions after \Cref{fundamental}.

We recall the standard notion of solution of \eqref{flow} first introduced in \cite{dipernalions}: we denote $\pazocal{B}(E,F)$ the space of bounded function from $E$ to $F$, $\log(L)(\mathbb{R}^d)$ the space
\[\left\{u\text{ measurable function}:\int_{\mathbb{R}^d} \log(1+|u(x)|)\,\dd x<\infty\right\},\]
$L^0(\mathbb{R}^d)$ the space of measurable functions with respect to Lebesgue measure $\mathcal{L}^d$. We are now able to introduce the ubiquitous notion of renormalized regular flow. The key insight is to consider the differential equation in with $\log(1+|\boldsymbol{X}(t,s,x)|)$ satisfies, where $\boldsymbol{X}$ solves \eqref{flow}.

\begin{definition}[Regular Lagrangian flow]\label{rlf}
	\textnormal{Let $s\in [0,T]$ and a vector field $\boldsymbol{b}$ satisfying \eqref{growth}. We say that  $\boldsymbol{X}:[s,T]\times[0,T]\times\mathbb{R}^d\rightarrow\mathbb{R}^d$ is a (renormalized) regular Lagrangian flow of \eqref{flow} starting at $t=s$ if for $D_T=\{(s,t)\in[0,T]^2: 0\leq s\leq t\leq T\}$, the regularity
		\[\boldsymbol{X}\in C(D_T;L^0_{\loc}(\mathbb{R}^d))\cap \pazocal{B}(D_T;\log(L)_{\loc}(\mathbb{R}^d)),\]
		holds, and $\boldsymbol{X}$ satisfies the renormalized flow equation in the distributional sense
		\begin{equation}\label{flowrenormalized}
			\begin{cases}
				&\partial_t[\beta(\boldsymbol{X}(t,s,x))]=\nabla\beta(\boldsymbol{X}(t,s,x))\cdot\boldsymbol{b}(t,\boldsymbol{X}(t,s,x))\\
				&\boldsymbol{X}(s,s,x)=x
			\end{cases}
		\end{equation}
		in $[s,T]\times\mathbb{R}^d$, where $\beta:\mathbb{R}^d\rightarrow \mathbb{R}$ is any function satisfying for some constant $C>0$ and almost every $x\in\mathbb{R}^d$
		\[|\beta(x)|\leq C(1+\log(1+|x|))\quad\text{and}\quad |\nabla\beta(x)|\leq C(1+|x|)^{-1}\quad \text{for all } x\in \mathbb{R}^d,\]
		and there exists a constant $L>0$ (called \textit{compressibility constant}) such that
		\[\int_{\mathbb{R}^d}\varphi(\boldsymbol{X}(t,s,y))\,\dd y\leq L\int_{\mathbb{R}^d}\varphi(y)\,\dd y\]
		for all measurable function $\varphi:\mathbb{R}^d\rightarrow \mathbb{R}_+$.
	}
\end{definition}

Our main result is the so called fundamental estimate for vector fields satisfying \eqref{growth} and \eqref{derivativeconvolution}. In particular, it ensure well-posedness of regular Lagrangian flows as presented in \cite{bouchutcrippa}.
\begin{theorem}[Fundamental estimate]\label{fundamental}
	Let $\boldsymbol{b},\,\bar{\boldsymbol{b}}\in L_{\loc}^1((0,T)\times\mathbb{R}^d))$ be vector fields satisfying \eqref{growth}. Moreover, assume that $\boldsymbol{b}$ with derivatives as in \eqref{derivativeconvolution}, with singular kernels $\Gamma_{ijk,R}$ as in \Cref{singularkernel} and summable functions $g^{ijk,R}$ in $(0,T)\times\mathbb{R}^d$ for some $q>1$ and for each index, and  $\boldsymbol{b}\in L^q_{\loc}((0,T)\times\mathbb{R}^d)$ . Finally, let $s,\tau\in [0,T]$ and $\boldsymbol{X},\bar{\boldsymbol{X}}$ be regular Lagrangian flows associated to $\boldsymbol{b},\bar{\boldsymbol{b}}$ starting at $t=s$ and with compressibility constants $L,\bar{L}$, respectively. Therefore for every $\gamma,r,\eta>0$, there exists $\lambda>0$ and a constant $C_{\gamma,r,\eta}>0$ depending on the $L^1((0,T)\times\mathbb{R}^d)$ norm of $g^{ijk,r}$, the $L^q((s,\tau)\times B_\lambda)$ norm of $\boldsymbol{b}$, the compressibility constants $L,\,\bar{L}$, and the $L^1((0,T)\times\mathbb{R}^d)$ norms of $\boldsymbol{b}^1,\,\bar{\boldsymbol{b}}^1$ and $L^1((0,T);L^\infty(\mathbb{R}^d))$ norms of $\boldsymbol{b}^2,\,\bar{\boldsymbol{b}}^2$, such that for any $t\in[s,T]$  it holds
	\begin{equation}\label{fundamentalestimate}
		\mathcal{L}^d\left(B_r\cap\{x\in\mathbb{R}^d:|\boldsymbol{X}(t,s,x)-\bar{\boldsymbol{X}}(t,s,x)|>\gamma\}\right)<C_{\gamma,r,\eta}\|\boldsymbol{b}-\bar{\boldsymbol{b}}\|_{L^1((s,\tau)\times B_\lambda)}+\eta. 
	\end{equation}
\end{theorem}

\Cref{fundamental} is a complete generalization of the result found in \cite{bouchutcrippa}, where the authors consider vector fields satisfying
\[\partial_j\boldsymbol{b}^i(t,x)=\sum^m_{k=1}\Gamma_{ijk,R}\ast g^{ijk,R}(t,x)\quad \text{for } (t,x)\in(0,T)\times B_R,\]
with $C^1$ kernels satisfying the appropriate point-wise decay
\[\left|\Gamma_{ijk,R}(x)\right|+|x|\left|\Gamma_{ijk,R}(x)\right|\leq C|x|^{-d}.\]
Hence, we are able to consider a much wider class of vector fields than \cite{bouchutcrippa}. Moreover, the recent result of Nguyen \cite{nguyen} does not contain (but neither is contained in) \Cref{fundamental}, as they have to assume a much stricter class of singular kernels and the structure of vector fields (not only for its derivatives). More precisely, they assume
\[\boldsymbol{b}^i(t,x)=\sum^m_{k=1}\Gamma_{ik,R}\ast g^{ik,R}(t,x)\quad \text{for } (t,x)\in(0,T)\times B_R,\]
where $g^{ik,R}\in L^1((0,T);\bv(\mathbb{R}^d))$ and $\Gamma_{ik,R}=|\cdot|^{-d}\Omega_{ik,R}$, $\Omega_{ik,R}\in \bv(\mathbb{S}^{d-1})$ being an average zero and zeroth order homogeneous function. Nevertheless, our result does not cover finite measures $g^{ijk,R}$, hence it is not an generalization of \cite{nguyen}.

As previously mentioned, the results in \cite[Sections 5 and 6]{bouchutcrippa} combined with \eqref{fundamentalestimate} implies the well-posedness of regular Lagrangian flows. Moreover, by the results in \cite[Section 7]{bouchutcrippa}, there exists a solution for the transport and continuity equations with vector fields $\boldsymbol{b}$.
\begin{proposition}\label{wellposedness}
	Let $\boldsymbol{b}$ as in \Cref{fundamental} and satisfying
		\[	\Div\boldsymbol{b}(t,x)\geq \alpha(t)\quad \text{in the weak sense for some } \alpha\in L^1([0,T]).\]
	Then for any $s\in[0,T]$, there exists an unique regular Lagrangian flow starting at $t=s$ as in \Cref{rlf} with semigroup property 
	\[\boldsymbol{X}(t,s,x)=\boldsymbol{X}(t,\tau,\boldsymbol{X}(\tau,s,x))\quad \text{for almost every } x\in\mathbb{R}^d \text{ and all } s\leq \tau\leq t\leq T.\] Moreover, there exists a solution of transport equation in the weak sense
	\[\begin{cases}
	\partial_t \beta(u)+\boldsymbol{b}\cdot\nabla \beta(u)=0\\
	u(0,\cdot)=u_0
	\end{cases}\]
	for any bounded $\beta\in C(\mathbb{R})$ and $u_0\in L^0(\mathbb{R}^d)$, provided that $\Div\boldsymbol{b}\in L^1_{\loc}((0,T)\times\mathbb{R}^d)$. The same holds for the continuity equation
	\[\begin{cases}
	\partial_t \beta(u)+\boldsymbol{b}\cdot\nabla \beta(u)+\beta'(u)u\Div\boldsymbol{b}=0\\
	u(0,\cdot)=u_0,
	\end{cases}\]
	for any bounded $\beta\in C^1(\mathbb{R})$ with $(1+|\cdot|)\beta'$ bounded and $u_0\in L^0(\mathbb{R}^d)$.
\end{proposition}

\textbf{Acknowledgments.} I would like to thank professor João Nariyoshi for the encouragement and insightful comments on the subject of this paper. This study was financed, in part, by the São Paulo Research Foundation (FAPESP), Brasil. Process Number 2024/21041-1. No datasets were generated or analysed during the current study. The authors have no competing interests to declare that are relevant to the content of this article.
\section{Preliminaries and proof of \texorpdfstring{\Cref{fundamental}}{}}
We first recall the definition of weak Lebesgue spaces $L^p_w(\mathbb{R}^d)$ for $p\in[1,\infty)$ as functions $u$ satisfying
\[\pseudo{u}_{L^p_w(\mathbb{R}^d)}\coloneqq\sup_{\lambda>0}\lambda\mathcal{L}^d(\{x\in\mathbb{R}^d: |u(x)|>\lambda\})^{\frac{1}{p}}<\infty.\]

We also recall the grand maximal operator majorly studied in \cite{bouchutcrippa}.
\begin{definition}\label{grandmaximal}
	\textnormal{ Let $u\in L^1_{\loc}(\mathbb{R}^d)$ and consider the family of functions $\{\rho^\nu\}_{\nu}\in L^1(\mathbb{R}^d)\cap L^\infty(\mathbb{R}^d)$ with $\supp\rho^\nu\subset B_1$. We define the grand maximal operator $\operatorname{M}_{\{\rho^\nu\}}$ as
		\[
		\operatorname{M}_{\{\rho^\nu\}}u(x)\coloneqq \sup_\nu\sup_{\epsilon>0}|\rho^\nu_\epsilon\ast u(x)|,
		\]
		where $\rho^\nu_\epsilon(x)=\epsilon^{-d}\rho^\nu(\epsilon^{-1}x)$.}
\end{definition}
Notice that
\[\operatorname{M}_{\{\rho^\nu\}}u(x)\leq\sup_\nu\|\rho^\nu\|_{L^\infty(\mathbb{R}^d)} \sup_{\epsilon>0}\frac{1}{\epsilon^d}\int_{B_\epsilon}|u(x-y)|\,\dd y\eqqcolon \sup_\nu\|\rho^\nu\|_{L^\infty(\mathbb{R}^d)}\operatorname{M}u(x),\]
and so the well-established $L^p(\mathbb{R}^d)$ to $L^p(\mathbb{R}^d)$ and $L^1(\mathbb{R}^d)$ to $L^1_w(\mathbb{R}^d)$ boundedness of maximal operator with $p\in(1,\infty)$ (for a classical reference, see \cite{stein}) also holds for the grand maximal operator. 
Also notice that one may consider $\rho(x)=\mathcal{L}^d(B_1)^{-1}\mathbbm{1}_{B_1}$ so that
\[\operatorname{M}_{\rho}u(x)= \sup_{\epsilon>0}\left|\frac{1}{\epsilon^d}\int_{B_\epsilon}u(x-y)\,\dd y\right|\eqqcolon \operatorname{M}u(x).\]

We state one of fundamental results in \cite[Lemma 2.4 and Remark 2.16]{nguyen}; although the result extends for finite measures, we shall state it for summable functions only for the sake of simplicity.
\begin{lemma}\label{limsupestimate}
	For $r>0$ and $q>1$, let $g\in L^q((0,T)\times B_r)$. Moreover, let $f\in L^1((0,T)\times\mathbb{R}^d)$ and $T$ be an operator bounded from $L^1(\mathbb{R}^d)$ to $L^1_w(\mathbb{R}^d)$, with limit superior
	\[
	\int_0^T\limsup_{\lambda\rightarrow \infty}\lambda\mathcal{L}^d(B_r\cap\{x\in \mathbb{R}^d:Tf(t,x)>\lambda\})\,\dd t=0,
	\]
	Then it holds
	\[\limsup_{\delta\rightarrow 0}\frac{1}{|\log\delta|}\int_0^T\int_{B_r}\min\left\{\frac{g(t,x)}{\delta}, Tf(t,x)\right\}\,\dd x\,\dd t=0.\]
\end{lemma}

Next, we state the control of sublevels
\begin{equation}\label{sublevel}
	G_\lambda\coloneqq \{x\in\mathbb{R}^d:|\boldsymbol{X}(t,s,x)|<\lambda \text{ for all } t\in[s,T]\}.
\end{equation}
first introduced in \cite[Proposition 3.2]{crippadelellis} and with a detailed proof in \cite[Lemma 5.5]{bouchutcrippa}.
\begin{lemma}\label{sublevelestimate}
	Let $s\in[0,T]$, $\boldsymbol{b}$ a vector field with growth assumption \eqref{growth}, and $\boldsymbol{X}$ a regular Lagrangian flow starting at $t=s$ with compressibility constant $L$ as in \Cref{rlf}. Then for any $r,\epsilon>0$, there exists $\lambda>0$ such that
	\[\mathcal{L}^d\left(B_r\setminus G_\lambda\right)<\epsilon.\]
\end{lemma}

We shall now prove \Cref{fundamental} under the assumption of \Cref{hormandertheorem} and \Cref{symmetric}; the proof is a direct adaptation of \cite[Proposition 5.9]{bouchutcrippa} and may be skipped by the familiar reader.
\begin{proof}[Proof of \Cref{fundamental}.] Considering for $0\leq s<t\leq T$ the function
	\[\Phi_\delta(t)\coloneqq \int_{G_\lambda\cap \bar{G}_\lambda\cap B_r}\log\left(1+\frac{|\boldsymbol{X}(t,s,x)-\bar{\boldsymbol{X}}(t,s,x)|}{\delta}\right)\,\dd x,\]
	where $r>0$ is arbitrary and $\delta,\lambda>0$ are to be chosen later, $\bar{G}_\lambda$ the sublevel with respect to $\bar{\boldsymbol{X}}$, we have as in \cite[proof of Proposition 5.9]{bouchutcrippa} that \Cref{sublevelestimate} implies for $\lambda$ large enough
	\begin{equation}\label{convergenceinmeasure}
		\begin{split}
		\mathcal{L}^d\left(B_r\cap\{x\in\mathbb{R}^d:|\boldsymbol{X}(t,s,x)-\bar{\boldsymbol{X}}(t,s,x)|>\gamma\}\right)\leq \frac{\Phi_\delta(t)}{\log\left(1+\frac{\gamma}{\delta}\right)}+\frac{\eta}{2}.
	\end{split}
	\end{equation}
	Moreover, by \Cref{symmetric} we may proceed as in \cite{bouchutcrippa} to obtain
	\[\begin{split}
		\Phi_\delta(t)\leq&\frac{\bar{L}}{\delta}\|\boldsymbol{b}-\bar{\boldsymbol{b}}\|_{L^1((s,t)\times B_\lambda)}\\&+\int_s^t\int_{G_\lambda\cap \bar{G}_\lambda\cap B_r}\frac{\boldsymbol{b}(\tau,\bar{\boldsymbol{X}}(\tau,s,x))-\boldsymbol{b}(\tau,\boldsymbol{X}(\tau,s,x))}{\delta+|\boldsymbol{X}(\tau,s,x)-\bar{\boldsymbol{X}}(\tau,s,x)|}\cdot\frac{\boldsymbol{X}(\tau,s,x)-\bar{\boldsymbol{X}}(\tau,s,x)}{|\boldsymbol{X}(\tau,s,x)-\bar{\boldsymbol{X}}(\tau,s,x)|}\,\dd x\,\dd \tau\\
		\leq&\frac{\bar{L}}{\delta}\|\boldsymbol{b}-\bar{\boldsymbol{b}}\|_{L^1((s,t)\times B_\lambda)}\\&+\int_s^t\int_{G_\lambda\cap \bar{G}_\lambda\cap B_r}\min\Bigg\{\frac{\left|\boldsymbol{b}(\tau,\bar{\boldsymbol{X}}(\tau,s,x))\right|+\left|\boldsymbol{b}(\tau,\boldsymbol{X}(\tau,s,x))\right|}{\delta},\\ &\quad\, U(\tau,\boldsymbol{X}(\tau,s,x))+U(\tau,\bar{\boldsymbol{X}}(\tau,s,x))\Bigg\}\,\dd x\,\dd \tau,
	\end{split}\]
	where we define
	\[\begin{split}
		U(\tau,z)\coloneqq& \operatorname{M}_{\Xi^\nu}(\sym(D \boldsymbol{b}(\tau,\cdot)))(z)+\operatorname{M}_{\Upsilon^\nu}(\sym(D \boldsymbol{b}(\tau,\cdot)))(z)+\operatorname{M}_{\Psi^\nu}(\sym(D \boldsymbol{b}(\tau,\cdot)))(z)\\&+\operatorname{M}_{\Omega^\nu}(\sym(D \boldsymbol{b}(\tau,\cdot)))(z)+\operatorname{M}_{\psi^\nu}(\Div\boldsymbol{b}(\tau,\cdot))(z)+\operatorname{M}_{\omega^\nu}(\Div\boldsymbol{b}(\tau,\cdot))(z).
	\end{split}\]
	Since we are assuming that $\sym(D\boldsymbol{b})$ (and consequently $\Div \boldsymbol{b}$) are of the form \eqref{derivativeconvolution}, one may apply \Cref{limsupestimate}, since by \Cref{hormandertheorem} we have that
	\begin{equation}\label{U}
		\left\|\pseudo{U}_{L^1_w(\mathbb{R}^d)}\right\|_{L^1(0,T)}\leq C_d(C_0+C_1)\sum_{i,j=1}^d\sum_{k=1}^m\|g^{ijk,r}\|_{L^1((0,T)\times\mathbb{R}^d)}
	\end{equation}
	and we are also assuming that $\boldsymbol{b}\in L^q_{\loc}((0,T)\times\mathbb{R}^d)$. Hence \eqref{convergenceinmeasure} gives that
	\[\mathcal{L}^d\left(B_r\cap\{x\in\mathbb{R}^d:|\boldsymbol{X}(t,s,x)-\bar{\boldsymbol{X}}(t,s,x)|>\gamma\}\right)\leq \eta+\frac{\|\boldsymbol{b}-\bar{\boldsymbol{b}}\|_{L^1((s,\tau)\times B_\lambda)}}{\delta_0\log\left(1+\frac{\gamma}{\delta_0}\right)},\]
	where $\delta_0>0$ depends on the list of constants in \Cref{fundamental}.
\end{proof}

Finally, we recall the ubiquitous Calderón-Zygmund decomposition; for a modern reference, see \cite[Section 5.3]{grafakos}.
\begin{theorem}[Calderón-Zygmund decomposition]\label{decomposition}
	Let $\alpha>0$ and $u\in L^1(\mathbb{R}^d)$. Then we may write $u=g+b$, where the ``good'' function $\|g\|_{L^1(\mathbb{R}^d)}\leq \|u\|_{L^1(\mathbb{R}^d)}$ and $\|g\|_{L^\infty(\mathbb{R}^d)}\leq 2^d\alpha$ and the ``bad'' function $b=\sum_{j}b_j$, with $\supp b_j\subset Q_j$ disjoint dyadic cubes $\{Q_j\}_j$. Moreover, $b_j$ has average zero,  $\sum_j\mathcal{L}^d(Q_j)\leq \alpha^{-1}\|u\|_{L^1(\mathbb{R}^d)}$, and $\|b_j\|_{L^1(Q_j)}\leq 2^{d+1}\alpha\mathcal{L}^d(Q_j)$. In particular, $\sum_j\|b_j\|_{L^1(Q_j)}\leq 2^{d+1}\|u\|_{L^1(\mathbb{R}^d)}$.
\end{theorem}
\section{Grand maximal operator and Hörmander condition}
We now shall prove the generalization of \cite[Theorem 3.3]{bouchutcrippa} for kernels as in \Cref{singularkernel}.
\begin{proposition}[Weak $L^1$ estimate]\label{hormandertheorem} Let $u\in L^1(\mathbb{R}^d)$, $\Gamma$ a singular kernel as in \Cref{singularkernel} and $\{\rho_\nu\}_{\nu}$ a family of functions as in \Cref{grandmaximal} such that
	\[C_2=\sup_\nu\|\rho^\nu\ast \Gamma\|_{L^\infty(\mathbb{R}^d)}<\infty.\]
Then there exists a constant $C$ depending only on the dimension $d$, $C_0$, $C_1$, and $C_2$ in \Cref{singularkernel} such that
\[\pseudo{\operatorname{M}_{\{\rho^\nu\}}\Gamma\ast u}_{L^1_w(\mathbb{R}^d)}\leq C\|u\|_{L^1(\mathbb{R}^d)};\]
moreover, if $\left\{\rho^\nu\right\}_\nu \in C_c^\infty(\mathbb{R}^d)$, then the aforementioned estimate holds for finite signed measures $\mu$, substituting $\|u\|_{L^1(\mathbb{R}^d)}$ with $\left|\mu\right|(\mathbb{R}^d)$.
\end{proposition}
Notice that $C_2$ being finite holds if the family $\{\rho^\nu\}\nu$ is in the fractional Hilbert space $H^{s}(\mathbb{R}^d)$ for $s>d/2$ by a straightforward application of Riemann-Lebesgue lemma.
\begin{proof} The proof follows standard techniques of singular kernel theory: write
	\[\begin{split}
		\rho^\nu\ast\Gamma\ast u(x)=&\frac{1}{\epsilon^d}\int_{\mathbb{R}^d}\int_{\mathbb{R}^d}\Gamma(z-y)\rho^\nu(\epsilon^{-1}y)u(x-z)\,\dd z\,\dd y\\=&\frac{1}{\epsilon^d}\int_{\mathbb{R}^d}\int_{\mathbb{R}^d}\mathbbm{1}_{|y-z|<2\epsilon}\Gamma(z-y)\rho^\nu(\epsilon^{-1}y)u(x-z)\,\dd z\,\dd y\\
		&+\frac{1}{\epsilon^d}\int_{\mathbb{R}^d}\int_{\mathbb{R}^d}\mathbbm{1}_{|y-z|\geq2\epsilon}\Gamma(z-y)\rho^\nu(\epsilon^{-1}y)u(x-z)\,\dd z\,\dd y\coloneqq \operatorname{I}_1(x)+\operatorname{I}_2(x).
	\end{split}\]
	Now, notice that the first integral is bounded as
	\[\operatorname{I}_1(x)\leq\frac{1}{\epsilon^d}\int_{|z|<3\epsilon}\left|\int_{|y|<2\epsilon}\Gamma(y)\rho^\nu(\epsilon^{-1}(z-y))\,\dd y\right||u(x-z)|\,\dd z,\]
	for the inner integral vanishes if $|z-y|> \epsilon$. Therefore by the characterization of singular kernels in \Cref{singularkernel}, which implies that $\mathcal{F}(\mathbbm{1}_{B_{2\epsilon}}\Gamma)\in L^\infty(\mathbb{R}^d)$ uniformly in $\epsilon$, combined with the finiteness of $C_2$ gives that
	\begin{equation}\label{I1}
		\operatorname{I}_1(x)\leq C_2\frac{1}{\epsilon^d}\int_{|z|<3\epsilon}|u(x-z)|\,\dd z\leq 3^dC_2 \operatorname{M}u(x).
	\end{equation}
	For the second integral, notice that
\begin{equation}\label{I2}
		\begin{split}
		\operatorname{I}_2(x)=&\frac{1}{\epsilon^d}\int_{\mathbb{R}^d}\int_{\mathbb{R}^d}\left[\mathbbm{1}_{|y-z|\geq2\epsilon}\Gamma(z-y)-\mathbbm{1}_{|z|\geq3\epsilon}\Gamma(z)\right]\rho^\nu(\epsilon^{-1}y)u(x-z)\,\dd z\,\dd y\\&+\|\rho^\nu\|_{L^1(\mathbb{R}^d)}\left(\mathbbm{1}_{|\cdot|>3\epsilon}\Gamma\right)\ast u(x)\eqqcolon \operatorname{I}_{21}(x)+\operatorname{I}_{22}(x).
	\end{split}
\end{equation}
	Notice that the classical theory of implies that \[\pseudo{\sup_{\epsilon>0}|\left(\mathbbm{1}_{|\cdot|>3\epsilon}\Gamma\right)\ast u|}_{L^1_w(\mathbb{R}^d)}\leq C\|u\|_{L^1(\mathbb{R}^d)},\]
	and so it suffices to obtain an analogous estimate for $\operatorname{I}_{21}$. For this purpose, notice that we only need to consider $|z|>\epsilon$, otherwise $\operatorname{I}_{21}$ vanishes. We thus write
	\[\begin{split}
		\{|z|>\epsilon\}=\left(\{|z|\geq3\epsilon,|z-y|\geq2\epsilon\}\right)\cup \left(\{|z|\geq3\epsilon,|z-y|<2\epsilon\}\right)\cup \left(\{\epsilon<|z|<3\epsilon,|z-y|\geq2\epsilon\}\right)\\\eqqcolon \Omega\cup \Omega'\cup \Omega'',
	\end{split}\]
	and notice that
	\begin{equation}\label{I21}
		\begin{split}
		\operatorname{I}_{21}(x)=& \iint_{\Omega}\left[\Gamma(z-y)-\Gamma(z)\right]\rho^\nu_\epsilon(y)u(x-z)\,\dd z\,\dd y-\iint_{\Omega'}\Gamma(z)\rho^\nu_\epsilon(y)u(x-z)\,\dd z\,\dd y\\&+\iint_{\Omega''}\Gamma(z-y)\rho^\nu_\epsilon(y)u(x-z)\,\dd z\,\dd y.
	\end{split}
	\end{equation}
	
	We now split the proof in step 1 for the integrals over $\Omega'$ and $\Omega''$; step 2 for the integral over $\Omega$; finally, in step 3 we conclude the proof and we repeat the argument in \cite[Theorem 3.3]{bouchutcrippa} for the measure case for the sake of completeness.
	
	\textbf{Step 1.} Notice that for $(y,z)\in\Omega'$, we have that $|y|>\epsilon$, thus the second integral in \eqref{I21} vanishes. Moreover, for $(y,z)\in \Omega''$ and $|y|<\epsilon$, we have that $2\epsilon\leq |y-z|<3\epsilon$. Therefore, since $\rho^\nu_\epsilon$ has support in $B_\epsilon$ we have that the integral vanishes if $|z|>4\epsilon$, and so
	\[\int_{2\epsilon\leq|y-z|<3\epsilon}|\rho^\nu_\epsilon(y)||\Gamma(z-y)|\,\dd y\leq \|\rho^\nu\|_{L^\infty(\mathbb{R}^d)}\left(\sup_{\delta>0}\int_{2\delta\leq|y|<3\delta}|\Gamma(y)|\,\dd y\right)\frac{1}{\epsilon^d}\mathbbm{1}_{B_{4\epsilon}}(z).\]
	Thus we conclude
	\[\iint_{\Omega''}\Gamma(z-y)\rho^\nu_\epsilon(y)u(x-z)\,\dd z\,\dd y\leq 3^d\|\rho^\nu\|_{L^\infty(\mathbb{R}^d)}\left(\sup_{\delta>0}\int_{2\delta\leq|y|<3\delta}|\Gamma(y)|\,\dd y\right) \operatorname{M}u(x).\]
	
	\textbf{Step 2.} For $(y,z)\in \Omega$, we shall write
	\[\operatorname{T}u(x)\coloneqq \sup_\nu\sup_{\epsilon>0}\left|\iint_{\Omega}[\Gamma(z-y)-\Gamma(z)]\rho^\nu_\epsilon(y)u(x-z)\,\dd z\,\dd y\right|.\]
	We claim that $\operatorname{T}$ is bounded operator from $L^1(\mathbb{R}^d)$ to $L_w^1(\mathbb{R}^d)$, with
	\[\pseudo{\operatorname{T}u}_{L^1_w(\mathbb{R}^d)}\leq C_d(C_0+C_1)\sup_{\nu}\left(\|\rho^\nu\|_{L^\infty(\mathbb{R}^d)}+\|\rho^\nu\|_{L^1(\mathbb{R}^d)}\right)\|u\|_{L^1(\mathbb{R}^d)}.\]
	For this purpose, we consider the Calderón-Zygmund decomposition \Cref{decomposition} for an finite combination of characteristic functions with disjoint dyadic intervals approximating $u$; notice that since such functions are dense in $L^1(\mathbb{R}^d)$, it suffices to prove the $L^1(\mathbb{R}^d)$ to $L^1_w(\mathbb{R}^d)$ boundedness of $\operatorname{T}$ for such functions. Moreover, we take the height on the decomposition as
	$\alpha'=\gamma_B\alpha$ with $\gamma_B=(2^{d+1}B)^{-1}$,
	where $B$ is the operator norm from $L^\infty(\mathbb{R}^d)$ to $L^\infty(\mathbb{R}^d)$ of $\operatorname{T}$; notice that $B<\infty$, as
	\begin{equation}\label{Bbound}
	\begin{split}
			\operatorname{T}u(x)\leq& \sup_\nu\|\rho^\nu\|_{L^1(\mathbb{R}^d)}\int_{|z-y|>2|y|}|\Gamma(z-y)-\Gamma(z)|\,\dd z\|u\|_{L^\infty(\mathbb{R}^d)}\\\leq& C_0\sup_\nu\|\rho^\nu\|_{L^1(\mathbb{R}^d)}\|u\|_{L^\infty(\mathbb{R}^d)},
	\end{split}
	\end{equation}
	where we have used \eqref{hormander}.
	Hence, by considering $Q_j^*$ the cube with parallel axis and same center as $Q^j$ and side length $2\sqrt{d}\;\ell(Q_j)$, $\ell(Q_j)$ being the side length of cube $Q_j$, we have that
	\begin{equation}\label{boundT}
		\begin{split}
		\mathcal{L}^d\left(\left\{x\in\mathbb{R}^d:\operatorname{T}u(x)\geq \alpha\right\}\right)&\leq \mathcal{L}^d\left(\left\{x\in\mathbb{R}^d:\operatorname{T}b(x)\geq \frac{\alpha}{2}\right\}\right)\\ &\leq \sum_i\mathcal{L}^d(Q_i^*)+\mathcal{L}^d\left(\left\{x\notin \cup_i Q_i^*:\operatorname{T}\sum_jb_j(x)\geq \frac{\alpha}{2}\right\}\right).
	\end{split}
	\end{equation}
	By \Cref{decomposition}, notice that the summation over $\{Q_i^*\}$ is bounded by $\alpha^{-1}\|u\|_{L^1(\mathbb{R}^d)}$. Therefore, we shall only consider $x\notin \cup_i Q_i^*$. For this purpose, denoting $z_j$ the center of $Q_j$, notice that
	\begin{equation}\label{Tbound}
		\begin{split}
		\operatorname{T}b_j(x)\leq& \sup_\nu\sup_{\epsilon>0}\left|\int_{Q_j\cap\{|x-z|>3\epsilon\}}b_j(z)\int_{|x-y-z|>2\epsilon}\rho^\nu_\epsilon(y)[\Gamma(x-y-z)-\Gamma(x-y-z_j)]\,\dd y\,\dd z\right|\\ &+ \sup_\nu\sup_{\epsilon>0}\left|\int_{Q_j\cap\{|x-z|>3\epsilon\}}b_j(z)\int_{|x-y-z|>2\epsilon}\rho^\nu_\epsilon(y)[\Gamma(x-y-z_j)-\Gamma(x-z_j)]\,\dd y\,\dd z\right|\\ &+ \sup_\nu\sup_{\epsilon>0}\left|\int_{Q_j\cap\{|x-z|>3\epsilon\}}b_j(z)\int_{|x-y-z|>2\epsilon}\rho^\nu_\epsilon(y)[\Gamma(x-z_j)-\Gamma(x-z)]\,\dd y\,\dd z\right|
	\end{split}
	\end{equation}
	Since $|x-z_j|\geq 2^{-1}\ell (Q_j^*)\geq \sqrt{d}\;\ell(Q_j)$ and $|z-z_j|\leq 2^{-1}\sqrt{d}\ell (Q_j)$, it holds that $|x-z_j|\geq 2|z-z_j|$. Therefore the third term summed over $j$ and then integrated over $\cup_i Q_i^*$ is bounded as
	\begin{equation}\label{thirdintegral}
		\begin{split}
		&\int_{\cup_i Q_i^*}\sum_j\sup_\nu\sup_{\epsilon>0}\left|\int_{Q_j\cap\{|x-z|>3\epsilon\}}b_j(z)\int_{|x-y-z|>2\epsilon}\rho^\nu_\epsilon(y)[\Gamma(x-z_j)-\Gamma(x-z)]\,\dd y\,\dd z\right|\,\dd x\\&\leq\sup_{\nu}\|\rho^\nu\|_{L^1(\mathbb{R}^d)}\sum_j\|b_j\|_{L^1(Q_j)}\int_{|x-z_j|\geq 2|z-z_j|}|\Gamma(x-z_j)-\Gamma(x-z)|\,\dd x\\&\leq 2^{d+1}\gamma_BC_0\sup_{\nu}\|\rho^\nu\|_{L^1(\mathbb{R}^d)}\|u\|_{L^1(\mathbb{R}^d)},
	\end{split}
	\end{equation}
	where we have used \eqref{hormander} and \Cref{decomposition}.
	
	For the second term, notice that $|x-z_j|>2|y|$ for all $x\notin \cup_i Q_i^*$ and $y\in B_\epsilon$, otherwise
	\[3\epsilon<|x-z|\leq 2|y|+|z-z_j|\leq 2\epsilon+\frac{\sqrt{d}\;\ell(Q_j)}{2}\implies \epsilon< \frac{\sqrt{d}\;\ell(Q_j)}{2};\]
	but $\sqrt{d}\;\ell(Q_j)\leq|x-z_j|\leq 2|y|\leq 2\epsilon$, a contradiction. Hence, the second term summed over $j$ and integrated over $\cup_i Q_i^*$ is estimated as
	\begin{equation}\label{secondintegral}
		\begin{split}
			&\int_{\cup_i Q_i^*}\sum_j\sup_\nu\sup_{\epsilon>0}\left|\int_{Q_j\cap\{|x-z|>3\epsilon\}}b_j(z)\int_{|x-y-z|>2\epsilon}\rho^\nu_\epsilon(y)[\Gamma(x-y-z_j)-\Gamma(x-z_j)]\,\dd y\,\dd z\right|\,\dd x\\&\leq\sup_{\nu}\|\rho^\nu\|_{L^1(\mathbb{R}^d)}\sum_j\|b_j\|_{L^1(Q_j)}\int_{|x-z_j|\geq 2|y|}|\Gamma(x-y-z_j)-\Gamma(x-z_j)|\,\dd x\\&\leq 2^{d+1}\gamma_BC_0\sup_{\nu}\|\rho^\nu\|_{L^1(\mathbb{R}^d)}\|u\|_{L^1(\mathbb{R}^d)}.
		\end{split}
	\end{equation}
	Therefore, we shall estimate the first term in \eqref{Tbound}. For this purpose, we may write it as
	\[\sup_\nu\sup_{\epsilon>0}\left|\int_{Q_j\cap\{|x-z|>3\epsilon\}}b_j(z)\int_{|y-z|>2\epsilon}\rho^\nu_\epsilon(x-y)[\Gamma(y-z)-\Gamma(y-z_j)]\,\dd y\,\dd z\right|.\]
	If $|y-z|\geq\sqrt{d}\;\ell(Q_j)$, then the result follows similarly as \eqref{secondintegral}, for $2|z-z_j|< \sqrt{d}\;\ell(Q_j)$, and so
	\begin{equation}\label{firstintegral1}
		\begin{split}
		&\sup_\nu\sup_{\epsilon>0}\left|\sum_j\int_{Q_j\cap\{|x-z|>3\epsilon\}}b_j(z)\int_{|y-z|>2\epsilon}\rho^\nu_\epsilon(x-y)[\Gamma(y-z)-\Gamma(y-z_j)]\,\dd y\,\dd z\right|\\&\leq \sup_\nu\|\rho^\nu\|_{L^\infty(\mathbb{R}^d)}\operatorname{M}\left(\sum_j\int_{Q_j}|b_j(z)| \mathbbm{1}_{|\cdot|>2|z-z_j|}|\Gamma(\cdot)-\Gamma(\cdot+(z-z_j))|\,\dd z\right)(x).
	\end{split}
	\end{equation}
	Hence, we may only consider $|y-z|< \sqrt{d}\;\ell(Q_j)$; notice that this implies that $2\epsilon<\sqrt{d}\;\ell(Q_j)$. Thus
	\[|y-z|\geq |x-z|-|x-y|\geq |x-z|-\epsilon\geq \frac{2}{3}|x-z|\geq \frac{1}{3}\sqrt{d}\;\ell(Q_j),\]
	where in the last inequality we have used that $2|z-z_j|<\sqrt{d}\;\ell(Q_j)<|x-z_j|$. Therefore, we have that
	\[\frac{1}{3}\sqrt{d}\;\ell(Q_j)\leq|y-z|\leq \sqrt{d}\;\ell(Q_j).\]
	Analogously, we have that
	\[|y-z_j|\leq |y-z|+|z-z_j|\leq \sqrt{d}\;\ell(Q_j)+\frac{1}{2}\sqrt{d}\;\ell(Q_j).\]
	Finally, we estimate from below $|y-z_j|$ as
	\[|y-z_j|\geq |x-z_j|-|y-x|\geq \sqrt{d}\;\ell(Q_j)-\epsilon\geq \frac{1}{2}\sqrt{d}\;\ell(Q_j).\]
	Summarizing, we have that
	\[\begin{cases}
		\frac{1}{3}\sqrt{d}\;\ell(Q_j)\leq|y-z|\leq \sqrt{d}\;\ell(Q_j);\\
		\frac{1}{2}\sqrt{d}\;\ell(Q_j)\leq|y-z_j|\leq \frac{3}{2}\sqrt{d}\;\ell(Q_j),
	\end{cases}\]
	and so we may write
	\begin{equation}\label{firstintegral2}
	\begin{split}
		&\sup_\nu\sup_{\epsilon>0}\left|\sum_j\int_{Q_j\cap\{|x-z|>3\epsilon\}}b_j(z)\int_{|y-z|>2\epsilon}\rho^\nu_\epsilon(x-y)[\Gamma(y-z)-\Gamma(y-z_j)]\,\dd y\,\dd z\right|\\&\leq \sup_\nu\|\rho^\nu\|_{L^\infty(\mathbb{R}^d)}\operatorname{M}\left(\sum_j\int_{Q_j}|b_j(z)| \mathbbm{1}_{\frac{1}{2}\sqrt{d}\;\ell(Q_j)\leq|\cdot-z_j|\leq \frac{3}{2}\sqrt{d}\;\ell(Q_j)}|\Gamma(\cdot-z_j)|\,\dd z\right)(x)\\
		&\quad+ \sup_\nu\|\rho^\nu\|_{L^\infty(\mathbb{R}^d)}\operatorname{M}\left(\sum_j\int_{Q_j}|b_j(z)| \mathbbm{1}_{\frac{1}{3}\sqrt{d}\;\ell(Q_j)\leq|\cdot-z|\leq \sqrt{d}\;\ell(Q_j)}|\Gamma(\cdot-z)|\,\dd z\right)(x).
	\end{split}
\end{equation}

 Combining \eqref{Tbound}, \eqref{thirdintegral}, \eqref{secondintegral}, \eqref{firstintegral1}, \eqref{firstintegral2} with \Cref{decomposition} and choosing $\gamma_B=(2^{d+1}B)^{-1}$ in the definition of $\alpha'=\gamma_B\alpha$ gives that
\[\mathcal{L}^d\left(\left\{x\in\mathbb{R}^d:\operatorname{T}u(x)\geq \alpha\right\}\right)\leq \frac{C_d}{\alpha}\left(C_0+C_1+B\right)\sup_{\nu}\left(\|\rho^\nu\|_{L^\infty(\mathbb{R}^d)}+\|\rho^\nu\|_{L^1(\mathbb{R}^d)}\right)\|u\|_{L^1(\mathbb{R}^d)}.\]
Since \eqref{Bbound} implies that $B<C_0\sup_{\nu}\|\rho^\nu\|_{L^1(\mathbb{R}^d)}$, the claim follows.

\textbf{Step 3.} Combining \eqref{I1} and \eqref{I2}, as well as \eqref{I21} with steps 1 and 2 implies the proposition with
\[C=C_d\left[C_2+(C_0+C_1)\sup_{\nu}\left(\|\rho^\nu\|_{L^\infty(\mathbb{R}^d)}+\|\rho^\nu\|_{L^1(\mathbb{R}^d)}\right)\right].\]
Finally, if $\{\rho^\nu\}_{\nu}\in C^\infty_c(\mathbb{R}^d)$ and we consider a finite signed measure $\mu$ in place of $u$, the result follows by density: consider $u_n\coloneqq\zeta_n\ast \mu$ an approximation of $\mu$, where $\zeta_n$ is a smooth mollifier. Notice that we may apply the first part of the proposition to $u_n$, so that
\[\pseudo{\operatorname{M}_{\{\rho^\nu\}}\Gamma\ast u_n}_{L^1_w(\mathbb{R}^d)}\leq C\|u_n\|_{L^1(\mathbb{R}^d)}\leq C|\mu|(\mathbb{R}^d).\]
Moreover, $\operatorname{M}_{\{\rho^\nu\}}\Gamma\ast u_n\rightarrow \operatorname{M}_{\{\rho^\nu\}}\Gamma\ast \mu$ in the sense of distributions as $n\rightarrow\infty$, which implies that for fixed $\nu$, $\epsilon$, and $x$, it holds $\rho^\nu_\epsilon\ast\Gamma\ast u_n(x)\rightarrow \rho^\nu_\epsilon\ast\Gamma\ast \mu(x)$ as $n\rightarrow\infty$. Hence for any $\lambda>0$
\[\mathbbm{1}_{\left\{\operatorname{M}_{\{\rho^\nu\}}\Gamma\ast \mu(x)>\lambda\right\}}\leq \liminf_{n\rightarrow \infty}\mathbbm{1}_{\left\{\operatorname{M}_{\{\rho^\nu\}}\Gamma\ast u_n(x)>\lambda\right\}},\]
and so the result follows by Fatou's Lemma.
\end{proof}
\section{Mean value inequality with symmetric derivative}
There have been a myriad of results proving an almost everywhere pointwise estimate 
\begin{equation}\label{previousmeanvalue}
	|u(x)-u(y)|\leq |x-y|\left(g(x)+g(y)\right)
\end{equation}
with $g$ generally depending on weak derivative $Du$, and we refer to \cite[Chapter II]{stein}, \cite[Proposition 4.2]{bouchutcrippa}, and \cite[Proposition 2.3]{nguyen}. More precisely, the classical result covered in \cite{stein} show that $g(x)=\operatorname{M}(D u)(x)$ is sufficient for \eqref{previousmeanvalue}. Moreover, the refinement first proved in \cite[Proposition 4.2]{bouchutcrippa} with $g(x)=\operatorname{M}_{\Upsilon^\nu}(D u)(x)$, where $\nu\in \mathbb{S}^{d-1}$ and $\Upsilon^\nu$ is a smooth vector field; we write the grand maximal function of a vector field $\boldsymbol{u}$ with respect to a family of vector fields $\{\Upsilon^\nu\}_\nu$ (with each component $\boldsymbol{u}_j$ as in \Cref{grandmaximal}) as the supremum of the inner product of convolutions, i.e. 
\[\operatorname{M}_{\Upsilon^\nu}(\boldsymbol{u})(x)=\sup_\nu\left|\sum_{j=1}^d\Upsilon^{\nu,j}\ast \boldsymbol{u}_j(x)\right|.\]

We shall use an analogous definition for matrices in the proof of \Cref{symmetric}: for a rank two tensor $U$ and a smooth family of rank two tensors $\{\Xi^\nu\}_\nu$ with each component $\Xi^\nu_{ij}$, we define
\begin{equation}\label{grandmaximalmatrix}
	\operatorname{M}_{\Xi^\nu}(U)(x)=\sup_\nu\left|\sum_{i,j=1}^d\Xi^{\nu,ij}\ast U_{ij}(x)\right|;
\end{equation}
recall that for matrices $A$ and $B$, the double inner product is defined as $A:B=\sum_{ij}A^{ij}B_{ij}$, and so \eqref{grandmaximalmatrix} is the supremum of the double inner product of convolutions.

Finally, in \cite[Proposition 2.3]{nguyen} they refine the previous result (in fact, they prove an general explicit equality) for $\rho>0$ small enough
\[g(x)= \operatorname{M}_{\Theta_1^{\xi,\rho}}\left(\frac{(x-y)}{|x-y|}\cdot D u\right)(x)+\rho\operatorname{M}_{\Theta_2^{\xi,\rho}}\left(D u\right)(x),\]
where $\xi\in \mathbb{S}^{d-1}$ and $\Theta_1^{\xi,\rho},\, \Theta_2^{\xi,\rho}$ are a family of functions and vector fields, respectively, with support in the cone \[\left\{x\in\mathbb{R}^d: \left|\frac{x}{|x|}-\xi\right|\leq \rho\right\}.\]

We shall prove a variation of \eqref{previousmeanvalue} akin to the one sided Lipschitz condition. More precisely, we shall prove that for vector fields $\boldsymbol{u}\in \bl(\mathbb{R}^d)$, it holds
\[(\boldsymbol{u}(x)-\boldsymbol{u}(y))\cdot(x-y)\leq |x-y|^2\left(g(x)+g(y)\right),\]
where now $g$ is a sum of various grand maximal operators (i.e., with different families of smooth functions) depending on the symmetric part of $D\boldsymbol{u}$; we shall denote it as $\sym(D\boldsymbol{u})$. Recall that for a matrix $A$, its symmetric part is defined as
\[\left[\sym(A)\right]_{ij}=\frac{1}{2}\left[A_{ij}+A_{ji}\right].\]

\begin{proposition}\label{symmetric}
		Let $\boldsymbol{u}\in L^1_{\loc}(\mathbb{R}^d)$ a vector field such that is weak symmetric  derivative $\sym(D\boldsymbol{u})$ such that $\operatorname{M}_{\rho^\nu}[\sym(D \boldsymbol{u})](x)<\infty$ for $x\in \mathbb{R}^d\setminus E$ and any family of functions $\{\rho^\nu\}_{\nu}$ as in \Cref{grandmaximal}, where $\mathcal{L}^d(E)=0$. Moreover, for $\alpha>0$, let $\zeta^\alpha$ be smooth normalized mollifier such that  $\supp \zeta\subset B_{\alpha}$, and assume that for a fixed $\alpha>0$, there exists a set $F_\alpha$ such that $\mathcal{L}^d(F_\alpha)=0$ and $\zeta^\alpha_\delta\ast \boldsymbol{u}_i(x)\rightarrow \boldsymbol{u}_i(x)$ for all $x\in \mathbb{R}^d\setminus F_\alpha$ and $i=1,\dots,d$, where $\zeta^\alpha_\delta(x)=\delta^{-d}\zeta^\alpha(\delta^{-1}x)$. Then for any normalized function $h\in C^\infty_c(\mathbb{R}^d)$  such that $\supp h\subset B_{2^{-1}}$ and for any $i,j=1,\dots,d$
		\begin{equation}\label{hcondition}
			w_j\partial_ih(w)-w_i\partial_jh(w)=0,
		\end{equation}
		 there exists a families of rank two tensors $\{\Xi^\nu,\Upsilon^\nu,\Psi^\nu, \Omega^\nu\}_{\nu\in\mathbb{S}^{d-1}}$ and $\psi^\nu=\tr\Psi^\nu$ $\omega^\nu=\tr\Omega^\nu$ such that for any $x,y\in  \mathbb{R}^d\setminus(E\cup F_\alpha)$, it holds
		\begin{equation}\label{symmetricestimate}
		\begin{split}
				(\boldsymbol{u}(x)-\boldsymbol{u}(y))\cdot(x-y)\leq |x-y|^2\Big(&\operatorname{M}_{\Xi^\nu}(\sym(D \boldsymbol{u}))(x)+\operatorname{M}_{\Xi^\nu}(\sym(D \boldsymbol{u}))(y)\\&+\operatorname{M}_{\Upsilon^\nu}(\sym(D \boldsymbol{u}))(x)+\operatorname{M}_{\Upsilon^\nu}(\sym(D \boldsymbol{u}))(y)\\
				&+\operatorname{M}_{\Psi^\nu}(\sym(D \boldsymbol{u}))(x)+\operatorname{M}_{\Psi^\nu}(\sym(D \boldsymbol{u}))(y)\\
				&+\operatorname{M}_{\Omega^\nu}(\sym(D \boldsymbol{u})(x)+\operatorname{M}_{\Omega^\nu}(\sym(D \boldsymbol{u})(y)\\
				&+\operatorname{M}_{\psi^\nu}(\Div\boldsymbol{u})(x)+\operatorname{M}_{\psi^\nu}(\Div\boldsymbol{u})(y)\\
				&+\operatorname{M}_{\omega^\nu}(\Div\boldsymbol{u})(x)+\operatorname{M}_{\omega^\nu}(\Div\boldsymbol{u})(y)\Big).
		\end{split}
		\end{equation}
\end{proposition}
We remark that if $h$ is a radial function, that is, if $h(w)=g(|w|)$ for some smooth $g\in C^\infty(\mathbb{R})$ such that $\supp g\in (-2^{-1},2^{-1})$, then \eqref{hcondition} holds.
\begin{proof}
	We proceed as in \cite[Proposition 4.2]{bouchutcrippa} and consider a normalized function $h\in C_c^\infty(\mathbb{R}^d)$ such that $\supp h\subset B_1$ and first prove the proposition for a smooth vector fields $\boldsymbol{u}$. For this purpose, denoting $h_r(z)=r^{-d}h(r^{-1}z)$ notice that
	\begin{equation}\label{splitsymmetric}
		\begin{split}
		(\boldsymbol{u}(x)-\boldsymbol{u}(y))\cdot(x-y)=&\int_{\mathbb{R}^d}h_r\left(z-\frac{x+y}{2}\right)(\boldsymbol{u}(x)-\boldsymbol{u}(z))\cdot(x-z)\,\dd z\\&-\int_{\mathbb{R}^d}h_r\left(z-\frac{x+y}{2}\right)(\boldsymbol{u}(y)-\boldsymbol{u}(z))\cdot(y-z)\,\dd z\\&+\int_{\mathbb{R}^d}h_r\left(z-\frac{x+y}{2}\right)(\boldsymbol{u}(x)-\boldsymbol{u}(z))\cdot(z-y)\,\dd z\\&-\int_{\mathbb{R}^d}h_r\left(z-\frac{x+y}{2}\right)(\boldsymbol{u}(y)-\boldsymbol{u}(z))\cdot(z-x)\,\dd z\\=&\operatorname{I}_1(x,y)+\operatorname{I}_2(x,y)+\operatorname{I}_3(x,y)+\operatorname{I}_4(x,y).
		\end{split}
	\end{equation}
	Notice that $\operatorname{I}_1(y,x)=-\operatorname{I}_2(x,y)$ and $\operatorname{I}_3(y,x)=-\operatorname{I}_4(x,y)$. We shall split the proof in step 1 for the estimate of $\operatorname{I}_1+\operatorname{I}_2$, step 2 for $\operatorname{I}_3+\operatorname{I}_4$, and finally step 3 for the nonsmooth case of $\boldsymbol{u}$.
	
	\textbf{Step 1.} Using the mean value theorem, a simple computation gives that
	\begin{equation}\label{I1forstep2}
		\begin{split}
		\operatorname{I}_1(x,y)=&\int_0^1\int_{\mathbb{R}^d}h_r\left(z-\frac{x+y}{2}\right)\sym(D\boldsymbol{u})(x+s(z-x)): (x-z)\otimes (x-z)\,\dd z\,\dd s\\=&\int_0^1\frac{1}{s^d}\int_{\mathbb{R}^d}h_r\left(\frac{x-y}{2}-\frac{w}{s}\right)\sym(D\boldsymbol{u})(x-w): \frac{w}{s}\otimes \frac{w}{s}\,\dd w\,\dd s\\=&r^2\int_0^1\frac{1}{(rs)^d}\int_{\mathbb{R}^d}h\left(\frac{x-y}{2r}-\frac{w}{rs}\right)\sym(D\boldsymbol{u})(x-w): \frac{w}{rs}\otimes \frac{w}{rs}\,\dd w\,\dd s,
	\end{split}
	\end{equation}
	and by the previous consideration, we have a similar expression for $\operatorname{I}_2(x,y)$. Taking $r=|x-y|$, we have that
	\begin{equation}\label{I1I2}
		\left|\operatorname{I}_1(x,y)+\operatorname{I}_2(x,y)\right|\leq |x-y|^2\left[\operatorname{M}_{\Xi^\nu}(\sym(D \boldsymbol{u}))(x)+\operatorname{M}_{\Xi^\nu}(\sym(D \boldsymbol{u}))(y)\right],
	\end{equation}
	where we have taken the supremum over $\epsilon=rs$ and for $\nu\in \mathbb{S}^{d-1}$, we define tensor
	\[\Xi^\nu(w)\coloneqq h\left(\frac{\nu}{2}-w\right)w\otimes w.\]
	
	\textbf{Step 2.} 	By the same symmetry as before, we now devote our proof for estimating $\operatorname{I}_3(x,y)$ and obtain an analogous one for $\operatorname{I}_4(x,y)$ by interchanging $x$ and $y$. We proceed as before and by mean value theorem, we have that
	\[\begin{split}
		\operatorname{I}_3(x,y)=&\sum_{i,j=1}^d\int_0^1\int_{\mathbb{R}^d}h_r\left(z-\frac{x+y}{2}\right)\partial_j\boldsymbol{u}_i(x+s(z-x))(x-z)_i(z-y)_j\,\dd z\,\dd s\\=&-\sum_{i,j=1}^d\int_0^1\frac{1}{s^d}\int_{\mathbb{R}^d}h_r\left(\frac{x-y}{2}-\frac{w}{s}\right)\partial_j\boldsymbol{u}_i(x-w)\left(\frac{w}{s}\right)_i\left(x-y-\frac{w}{s}\right)_j\,\dd z\,\dd s.
	\end{split}\]
	We now write
	\[\begin{split}
		-\left(\frac{w}{s}\right)_i\left(x-y-\frac{w}{s}\right)_j=&\left[\left(\frac{w}{s}\right)_i\left(\frac{x-y}{2}\right)_j+\left(\frac{x-y}{2}\right)_i\left(\frac{w}{s}\right)_j\right]+\left(\frac{w}{s}\right)_i\left(\frac{w}{s}\right)_j\\&+\left[\left(\frac{x-y}{2}\right)_i\left(\frac{w}{s}\right)_j-\left(\frac{w}{s}\right)_i\left(\frac{x-y}{2}\right)_j\right]
	\end{split}\]
	and so we split $\operatorname{I}_3(x,y)$ as
	\[\begin{split}
			&\operatorname{I}_3(x,y)\\&=-\sum_{i,j=1}^d\int_0^1\frac{1}{s^d}\int_{\mathbb{R}^d}h_r\left(\frac{x-y}{2}-\frac{w}{s}\right)\left[\left(\frac{w}{s}\right)_i\left(\frac{x-y}{2}\right)_j+\left(\frac{x-y}{2}\right)_i\left(\frac{w}{s}\right)_j\right]\,\dd z\,\dd s\\&\quad+\sum_{i,j=1}^d\int_0^1\frac{1}{s^d}\int_{\mathbb{R}^d}h_r\left(\frac{x-y}{2}-\frac{w}{s}\right)\partial_j\boldsymbol{u}_i(x-w)\left(\frac{w}{s}\right)_i\left(\frac{w}{s}\right)_j\,\dd z\,\dd s\\&\quad+\sum_{i,j=1}^d\int_0^1\frac{1}{s^d}\int_{\mathbb{R}^d}h_r\left(\frac{x-y}{2}-\frac{w}{s}\right)\partial_j\boldsymbol{u}_i(x-w)\left[\left(\frac{x-y}{2}\right)_i\left(\frac{w}{s}\right)_j-\left(\frac{w}{s}\right)_i\left(\frac{x-y}{2}\right)_j\right]\,\dd z\,\dd s
	\end{split}\]
	Notice that the second term is the same as in step 1. Moreover, proceeding as in step 1, the second one can be written only in terms of $\sym(D\boldsymbol{u})$, for
	\[\begin{split}
		&\left|\sum_{i,j=1}^d\int_0^1\frac{1}{s^d}\int_{\mathbb{R}^d}h_r\left(\frac{x-y}{2}-\frac{w}{s}\right)\partial_j\boldsymbol{u}_i(x-w)\left[\left(\frac{w}{s}\right)_i\left(\frac{x-y}{2}\right)_j+\left(\frac{x-y}{2}\right)_i\left(\frac{w}{s}\right)_j\right]\,\dd z\,\dd s\right|\\&\leq |x-y|^2\operatorname{M}_{\Upsilon^\nu}(\sym(D\boldsymbol{u}))(x),
	\end{split}\]
	where we define for $\nu\in\mathbb{S}^{d-1}$ the rank two tensor
	\[\Upsilon^\nu(w)=\frac{1}{2}h\left(\frac{\nu}{2}-w\right)\left[w\otimes \nu+\nu\otimes w\right].\]
	recall that $r=|x-y|$. Hence, we conclude that
	\begin{equation}\label{I3}
		\left|\operatorname{I}_3(x,y)\right|\leq |x-y|^2\left(\operatorname{M}_{\Xi^\nu}\sym{D\boldsymbol{u}}(x)+\operatorname{M}_{\Upsilon^\nu}\sym{D\boldsymbol{u}}(x)\right)+\left|\operatorname{I}_5(x,y)\right|,
	\end{equation}
	where we define
	\[\operatorname{I}_5(x,y)\coloneqq\sum_{i,j=1}^d\int_0^1\frac{1}{s^d}\int_{\mathbb{R}^d}h_r\left(\frac{x-y}{2}-\frac{w}{s}\right)\partial_j\boldsymbol{u}_i(x-w)\left[\left(\frac{x-y}{2}\right)_i\left(\frac{w}{s}\right)_j-\left(\frac{w}{s}\right)_i\left(\frac{x-y}{2}\right)_j\right]\,\dd z\,\dd s\]
	and we shall devote the rest of this step into estimating the aforementioned term. Now, notice that $2w_i=\partial_i|w|^2$, and so integrating by parts and using \eqref{hcondition}, we have that
	\[\begin{split}
		&\operatorname{I}_5(x,y)\\&=\sum_{i,j=1}^d\int_0^1\frac{1}{s^d}\int_{\mathbb{R}^d}\left[\left(\frac{w}{s}\right)_i\partial_jh_r\left(\frac{x-y}{2}-\frac{w}{s}\right)-\left(\frac{w}{s}\right)_j\partial_ih_r\left(\frac{x-y}{2}-\frac{w}{s}\right)\right]\frac{|w|^2}{2rs^2}\partial_j\boldsymbol{u}_i(x-w)\,\dd z\,\dd s\\&\quad+\sum_{i,j=1}^d\int_0^1\frac{1}{s^d}\int_{\mathbb{R}^d}h_r\left(\frac{x-y}{2}-\frac{w}{s}\right)\left(\frac{x-y}{2}\right)_i\partial_{jj}\boldsymbol{u}_i(x-w)\frac{|w|^2}{2s}\,\dd z\,\dd s\\&\quad-\sum_{i,j=1}^d\int_0^1\frac{1}{s^d}\int_{\mathbb{R}^d}h_r\left(\frac{x-y}{2}-\frac{w}{s}\right)\left(\frac{x-y}{2}\right)_j\partial_{ji}\boldsymbol{u}_i(x-w)\frac{|w|^2}{2s}\,\dd z\,\dd s.
	\end{split}\]
	Recalling the definition of divergence of $U$ rank two tensor $(\Div U)_i=\sum_{j=1}^d\partial_j U_{ij}$, the vector calculus identity
	\begin{equation}\label{vectorcalculus}
		\Delta \boldsymbol{u}_i=2[\Div \sym(D\boldsymbol{u})]_i-\partial_i\Div\boldsymbol{u},
	\end{equation}
	defining for $\nu\in\mathbb{S}^{d-1}$ the smooth rank two tensor
	\[\Psi^\nu(w)=h\left(\frac{\nu}{2}-w\right)w\otimes\nu-\frac{|w|^2}{2}\nabla h\left(\frac{\nu}{2}-w\right)\otimes \nu,\]
	and denoting $\psi^\nu=\tr(\Psi^\nu)$,
	the last two terms can be estimated as
	\begin{equation}\label{I5estimate}
		\begin{split}
			&\left|\sum_{i,j=1}^d\int_0^1\frac{1}{s^d}\int_{\mathbb{R}^d}h_r\left(\frac{x-y}{2}-\frac{w}{s}\right)\left(\frac{x-y}{2}\right)_i\partial_{jj}\boldsymbol{u}_i(x-w)\frac{|w|^2}{2s}\,\dd z\,\dd s\right|\\&+\left|\sum_{i,j=1}^d\int_0^1\frac{1}{s^d}\int_{\mathbb{R}^d}h_r\left(\frac{x-y}{2}-\frac{w}{s}\right)\left(\frac{x-y}{2}\right)_j\partial_{ji}\boldsymbol{u}_i(x-w)\frac{|w|^2}{2s}\,\dd z\,\dd s\right|\\
			&\leq|x-y|^2\left[\operatorname{M}_{\Psi^\nu}(\sym(D\boldsymbol{u}))(x)+\operatorname{M}_{\psi^\nu}(\Div\boldsymbol{u})\right].
		\end{split}
	\end{equation}
	
	Hence, we shall complete step 2 if we obtain the desired estimate for
	\[\begin{split}
		&\operatorname{I}_6(x,y)\\&\coloneqq\sum_{i,j=1}^d\int_0^1\frac{1}{s^d}\int_{\mathbb{R}^d}\left[\left(\frac{w}{s}\right)_i\partial_jh_r\left(\frac{x-y}{2}-\frac{w}{s}\right)-\left(\frac{w}{s}\right)_j\partial_ih_r\left(\frac{x-y}{2}-\frac{w}{s}\right)\right]\frac{|w|^2}{2rs^2}\partial_j\boldsymbol{u}_i(x-w)\,\dd z\,\dd s.
	\end{split}\]
	For this purpose, using once again $2w_i=\partial_i|w|^2$, we have that
	\[\begin{split}
		&\sum_{i,j=1}^d\int_0^1\frac{1}{s^d}\int_{\mathbb{R}^d}\left(\frac{w}{s}\right)_i\partial_jh_r\left(\frac{x-y}{2}-\frac{w}{s}\right)\frac{|w|^2}{2rs^2}\partial_j\boldsymbol{u}_i(x-w)\,\dd z\,\dd s\\&=-\sum_{i,j=1}^d\int_0^1\frac{1}{s^d}\int_{\mathbb{R}^d}\frac{|w|^2}{2rs^2}\partial_jh_r\left(\frac{x-y}{2}-\frac{w}{s}\right)\left(\frac{w}{s}\right)_i\partial_j\boldsymbol{u}_i(x-w)\,\dd z\,\dd s\\&\quad+\sum_{i,j=1}^d\int_0^1\frac{1}{s^d}\int_{\mathbb{R}^d}\partial_{ij}h_r\left(\frac{x-y}{2}-\frac{w}{s}\right)\frac{|w|^4}{4r^2s^4}\partial_j\boldsymbol{u}_i(x-w)\,\dd z\,\dd s\\&\quad+\sum_{i,j=1}^d\int_0^1\frac{1}{s^d}\int_{\mathbb{R}^d}\partial_{j}h_r\left(\frac{x-y}{2}-\frac{w}{s}\right)\frac{|w|^4}{4rs^3}\partial_{ji}\boldsymbol{u}_i(x-w)\,\dd z\,\dd s,
	\end{split}\]
	and so we may write it as
	\begin{equation}\label{integrationbyparts1}
		\begin{split}
			&\sum_{i,j=1}^d\int_0^1\frac{1}{s^d}\int_{\mathbb{R}^d}\left(\frac{w}{s}\right)_i\partial_jh_r\left(\frac{x-y}{2}-\frac{w}{s}\right)\frac{|w|^2}{2rs^2}\partial_j\boldsymbol{u}_i(x-w)\,\dd z\,\dd s\\&=\frac{1}{2}\sum_{i,j=1}^d\int_0^1\frac{1}{s^d}\int_{\mathbb{R}^d}\partial_{ij}h_r\left(\frac{x-y}{2}-\frac{w}{s}\right)\frac{|w|^4}{4r^2s^4}\partial_j\boldsymbol{u}_i(x-w)\,\dd z\,\dd s\\&\quad+\frac{1}{2}\sum_{i,j=1}^d\int_0^1\frac{1}{s^d}\int_{\mathbb{R}^d}\partial_{j}h_r\left(\frac{x-y}{2}-\frac{w}{s}\right)\frac{|w|^4}{4rs^3}\partial_{ji}\boldsymbol{u}_i(x-w)\,\dd z\,\dd s.
		\end{split}
	\end{equation}
	Likewise, we may write
		\begin{equation}\label{integrationbyparts2}
		\begin{split}
			&\sum_{i,j=1}^d\int_0^1\frac{1}{s^d}\int_{\mathbb{R}^d}\left(\frac{w}{s}\right)_j\partial_ih_r\left(\frac{x-y}{2}-\frac{w}{s}\right)\frac{|w|^2}{2rs^2}\partial_j\boldsymbol{u}_i(x-w)\,\dd z\,\dd s\\&=\frac{1}{2}\sum_{i,j=1}^d\int_0^1\frac{1}{s^d}\int_{\mathbb{R}^d}\partial_{ij}h_r\left(\frac{x-y}{2}-\frac{w}{s}\right)\frac{|w|^4}{4r^2s^4}\partial_j\boldsymbol{u}_i(x-w)\,\dd z\,\dd s\\&\quad+\frac{1}{2}\sum_{i,j=1}^d\int_0^1\frac{1}{s^d}\int_{\mathbb{R}^d}\partial_{i}h_r\left(\frac{x-y}{2}-\frac{w}{s}\right)\frac{|w|^4}{4rs^3}\partial_{jj}\boldsymbol{u}_i(x-w)\,\dd z\,\dd s.
		\end{split}
	\end{equation}
	Combining \eqref{integrationbyparts1} and \eqref{integrationbyparts2} and applying \eqref{vectorcalculus}, we have that
	\begin{equation}\label{I5estimate1}
		\left|\operatorname{I}_6(x,y)\right|\leq |x-y|^2\left[\operatorname{M}_{\Omega^\nu}(\sym(D\boldsymbol{u}))(x)+\operatorname{M}_{\omega^\nu}(\Div\boldsymbol{u})(x)\right],
	\end{equation}
	where we define for $\nu\in\mathbb{S}^{d-1}$ the rank two tensor
	\[\Omega^\nu(w)\coloneqq \frac{|w|^4}{4}D^2h\left(\frac{\nu}{2}-w\right)+|w|^2\nabla h\left(\frac{\nu}{2}-w\right)\otimes w\]
	and denote $\omega^\nu\coloneqq \tr(\Omega^\nu)$. Combining steps 1 and 2, namely \eqref{splitsymmetric}, \eqref{I1I2}, \eqref{I3}, \eqref{I5estimate}, and \eqref{I5estimate1}, we obtain \eqref{symmetricestimate} for smooth vector fields $\boldsymbol{u}$.
	 
	\textbf{Step 3.} For the general case, we proceed by approximation: choose $\alpha>0$ in the hypothesis of \Cref{symmetric} so that $\supp h\subset B_{2^{-1}-\alpha}$ and let $\boldsymbol{u}_\delta=\zeta_\delta\ast \boldsymbol{u}$ (meaning that $(\boldsymbol{u}_\delta)_i=\zeta_\delta\ast \boldsymbol{u}_i$ for each index $i=1,\dots, d$). Since $\boldsymbol{u}_\delta$ is smooth, we apply \eqref{splitsymmetric} to it and by the associativity of convolution, we have an approximate version of \eqref{I1forstep2}:
	\begin{equation}\label{approximation}
		\operatorname{I}_1(x,y)=|x-y|^2\int_0^1\frac{1}{(rs)^d}\Xi^{\frac{x-y}{|x-y|},\delta}\left((sr)^{-1}\cdot\right)\ast\sym(D\boldsymbol{u})(x)\,\dd s,
	\end{equation}
	where the convolution is once again meant component-wise and $\Xi^{\nu,\delta}(w)\coloneqq\zeta_\delta\ast \Xi^\nu(w)$.
	
	Now, we may write as in \cite[Step 2 of Proposition 4.2]{bouchutcrippa}
	\[(\Xi^{\nu,\delta})_\epsilon(w)=(\zeta_\delta\ast\Xi^{\nu})_\epsilon(w)=\begin{cases}
		(\zeta\ast \Xi^\nu_{\delta^{-1}\epsilon})_\delta(w)\quad&\text{if } \epsilon\leq \delta;\\
		(\zeta_{\epsilon^{-1}\delta}\ast \Xi^\nu)_\epsilon(w) \quad&\text{if } \epsilon>\delta.
	\end{cases}\]
Hence, considering the smooth family
	\[\Pi_{\Xi}^{\nu,\epsilon,\delta}(w)=\begin{cases}
		(\zeta\ast \Xi^\nu_{\delta^{-1}\epsilon})_\delta(w)\quad&\text{if } \epsilon\leq \delta;\\
		(\zeta_{\epsilon^{-1}\delta}\ast \Xi^\nu)_\epsilon(w) \quad&\text{if } \epsilon>\delta,
	\end{cases}\]
	we have that $\supp \Pi_\Xi^{\nu,\epsilon,\delta}\subset B_1$, and so the family $\{\Pi_\Xi^{\nu'}\}_{\nu'}$ with $\nu'=(\nu,\epsilon,\delta)$ satisfies the hypothesis on \Cref{grandmaximal}, hence there exists a measure zero set $E_\Xi$ such that $\operatorname{M}_{\Pi_{\Xi}^{\nu'}}(\sym(D\boldsymbol{u}))(x)<\infty$ for all $x\in \mathbb{R}^d\setminus E_\Xi$. Moreover, by hypothesis there exists a measure zero set $F_\alpha$ such that $\boldsymbol{u}_\delta(x)\rightarrow \boldsymbol{u}(x)$ for all $x\in\mathbb{R}^d\setminus F_\alpha$, and so by dominated convergence theorem in \eqref{approximation} that we may pass to the limit as $\delta\rightarrow 0$
	for all $x\in\mathbb{R}^d\setminus(E_\Xi\cup F_\alpha)$ and \eqref{I1I2} still follows.
	
	The same procedure can be done analogously for the other kernels 	\[\Upsilon^{\nu,\delta}(w)\coloneqq\zeta_\delta\ast \Upsilon^\nu(w),\quad \Psi^{\nu,\delta}(w)\coloneqq\zeta_\delta\ast \Psi^\nu(w),\quad \Omega^{\nu,\delta}(w)\coloneqq\zeta_\delta\ast \Omega^\nu(w),\]
	as well as $\psi^{\nu,\delta}\coloneqq \tr\Psi^{\nu,\delta}$ and $\omega^{\nu,\delta}\coloneqq \tr\Omega^{\nu,\delta}$. Thus, there exists a set\[E\coloneqq E_\Xi\cup E_\Upsilon\cup E_\Psi\cup E_\Omega \cup E_\psi \cup E_\omega\]
	in which \eqref{symmetricestimate} holds for all $x,y\in \mathbb{R}^d\setminus (E\cup F_\alpha)$.
\end{proof}
\bibliographystyle{plain}
\bibliography{bib}
\end{document}